**Note on Integer Factoring Algorithms II**
**N. A. Carella, February 2006**

**Abstract:** This note introduces a new class of integer factoring algorithms. Two versions of this method will be described, deterministic and probabilistic. These algorithms are practical, and can factor large classes of balanced integers $N = pq$, $p < q < 2p$ in superpolynomial time. Further, an extension of the Fermat factoring method is proposed.

# 1 Introduction

This note introduces a new class of integer factorization algorithms. Two versions of this method will be described, one deterministic and one probabilistic. These techniques are practical, and can factor large classes of weak integers indistinguishable from random balanced integers $N = pq$, $p < q < 2p$, in superpolynomial time. These techniques exploit the additive structures (binary expansions, et cetera) of the integer of $p$, $q$, $p + q$, and $q - p$. These are the counterparts of the $p - 1$ method, $p + 1$ method and the elliptic curve method based on the multiplicative structures of the integers $p - 1$, $p + 1$, and $p + 1 - a_q$, where $| a_p | < 2p^{1/2}$, respectively. Further, an extension of the Fermat factoring method to all odd integers is proposed. It has deterministic time complexity of $O(N^{1/4})$ arithmetic operations, and it appears to be new in the literature.

Section 2 presents an extension of the Fermat factoring method to all odd composite integers $N$, see Theorem 4. Section 3 describes the *sparse difference method* (Theorem 7) and related ideas. Section 4 describes the *sparse exponent method*, see Theorem 10. A summary of the new classes of weak balanced integers appears in Section 5. Sections 6 and 7 are optional and are included as references. Some definitions, notations, and a few well known factorization methods are recorded in Section 6. An introduction to the foundation of the new integer factorization methods, and some useful number theoretical results are recorded in Section 7.

# 2 Extended Fermat Method

This section addresses the factorization of integers $N = pq$ as a function of either the factors sum $p + q$ or difference $q - p$. The factors sum $p + q$ and difference $q - p$ are equivalent approach linked by the identity $4pq = (p + q)^2 - (q - p)^2$. The earliest factoring algorithm based on this identity appears to be the Fermat method, also known as the difference of squares method, [RL, p. 147]. Magically, this factoring algorithm splits an odd integer $N$ in a few arithmetic steps whenever there is a pair of factors $p$ and $q$ close to the center $\sqrt{N}$ of $N$ or equivalently the factors difference $q - p \leq O(N^{1/4})$ is small.

For every odd integer $N$ the equation $N = x^2 - y^2$ has a finite number of solutions $(x, y)$ of the form $x = (p + q)/2$, $y = (p - q)/2$ with $p, q \mid N$. The factors of an arbitrary odd integer $N$, which can be composites or primes, vary from $p = q = N^{1/2}$, to $p = N/3$, $q = 3$, and $p = N$, $q = 1$. Equality occurs if and only if $N$ is a square.

The solution $N = (2^{e-1}m + 1)^2 - (2^{e-1}m)^2$ or $(2^{e-1}m)^2 - (2^{e-1}m - 1)^2$ obtained from $N = 2^e m \pm 1$, $m$ odd, $e \geq 1$, is viewed as the trivial solution. Prime numbers $N > 2$ have only trivial solutions. Now if $N$ is not prime, then there is a nontrivial solution such that $x \geq N^{1/2}$ is an integer in the sequence of numbers

$$x_0 = \sqrt{N}, \ x_1 = \sqrt{N} + 1, x_2 = \sqrt{N} + 2, \ ..., \ x_n = (N + 9)/6. \tag{1}$$

An inspection of the sequence (1) shows that the number of steps required to find a solution is

$$p + N/p - 2\sqrt{N} = (p - \sqrt{N})^2 / p. \tag{2}$$

Clearly, if there is a factor $p$ sufficiently close to $N^{1/2}$ the procedure is successful, and runs in deterministic polynomial time (even constant time). Otherwise,

$$p < N^{1/2} - N^{1/4+\alpha} \quad \text{and} \quad N^{1/2} + N^{1/4+\alpha} < q \tag{3}$$

for some $\alpha > 0$, and the algorithm runs in exponential time $O(N^{2\alpha})$.

***Theorem* 1.** Let $X > 0$ be a sufficiently large number and let $N \leq X$. Define the variable $Z(N) = \#\{ (x, y) \in \mathbb{Z} \times \mathbb{Z} : x^2 - y^2 = N \}$. Then the followings hold.
(i) The decision problem $Z(N) = 0$, or 1 or $> 1$ has deterministic time complexity.
(ii) The random variable $Z(N) \leq O(\log N)$ has a Gaussian distribution.

Proof: The first statement follows from the Agrawal-Kayal-Saxena primality test, and the second statement follows from the Erdos-Kac Theorem. ∎

The decision problem (i) for an arbitrary quadratic equation $ax^2 + bxy + cy^2 = N$ seems to have exponential time complexity [L1,2], and for cubic equations it is almost completely unknown.

The set of composite integers factorable by the Fermat method and the associated counting function are given by

$$\mathcal{F} = \{ N = pq : q - p = O(N^{1/4}) \} \ \text{ and } \ \mathcal{F}(X) = \#\{ N = pq \leq X : q - p = O(N^{1/4}) \} \tag{4}$$

respectively. Similarly, the set of balanced integers and the associated counting function are denoted by $\mathcal{B} = \{ N = pq : p < q < 2p \}$ and $\mathcal{B}(X) = \#\{ N = pq \leq X : p < q < 2p \}$ respectively.

***Theorem* 2.** $\mathcal{F}(X) \leq X^{3/4}/(\log X)^2$.

Proof: Since the pair $x = p + d$, $y = d$ is a solution of $N = x^2 - y^2$, this amounts to a count of the prime pairs $p \leq X^{1/2}$ and $q = p + 2d$ with $2d \leq X^{1/4}$. A rough estimate is given by

$$\mathcal{F}(X) = \pi(X^{1/2})[\pi(X^{1/2} + X^{1/4}) - \pi(X^{1/2})] \leq \left(\frac{c_0 X^{1/2}}{\log X^{1/2}}\right)\left(\frac{c_1 X^{1/4}}{\log X^{1/4}}\right) = \frac{c_2 X^{3/4}}{(\log X)^2}, \qquad (5)$$

where $\pi(X^{1/2}) = c_0 X^{1/4}/(\log X^{1/2})$, and $\pi(X^{1/2} + X^{1/4}) - \pi(X^{1/2}) \leq c_1 X^{1/4}/(\log X^{1/4})$ is the maximal number of primes in the interval $[X^{1/2}, X^{1/2} + X^{1/4}]$, and $c_0, c_1, c_2$ are constants. ∎

Since $\mathcal{B}(X) = c_3 X/(\log X)^2$, it follows that the subset $\mathcal{F}$ has zero density in the set $\mathcal{B}$. A more refined and general approach by sieve methods yields about the same result, see Theorem 19.

***Theorem* 3.** Let $N$ be an odd integer, and let $p = N^{1/2} + aN^{1/4} + b$ be a factor of $N$, where $0 \leq |a|, |b| < N^{1/4}$. Then the complexity of the Fermat method is as follows:
(i) Given the approximation $N^{1/2} + aN^{1/4}$ of $p$, the algorithm runs in deterministic time $O(ab/N^{1/4})$.
(ii) If $|a| \leq O(N^{\varepsilon})$ and $|b| \leq O(N^{1/4-\varepsilon})$, $\varepsilon > 0$, then the algorithm runs in deterministic time $O(N^{\varepsilon})$.
(iii) If $a \leq O(\log(N)^c)$, $c > 0$ constant, then the algorithm runs in deterministic polynomial time.

Proof (i): Let $N^{1/2} + aN^{1/4}$ be an approximation of a factor $p$ of $N$. Then the number of steps required to find an integral solution of $4N = x^2 - y^2$ is given by

$$p + \frac{N}{p} - \left(N^{1/2} + aN^{1/4} + \frac{N}{N^{1/2} + aN^{1/4}}\right) = \frac{(2aN^{3/4} + (a^2 + b)N^{1/2} + aN^{1/4})}{(N^{1/2} + aN^{1/4})(N^{1/2} + aN^{1/4} + b)} b. \qquad (6)$$

Hence, for sufficiently large $N$ and $a \neq 0$, the number of steps $\leq |4ab/N^{1/4}| \leq 4N^{1/4}$.
Statements (ii) and (iii) follow from (i). In particular, $a = 0$ reduces to the classical case described in (1) and (2). ∎

The previous attempts in the literature to extend the Fermat method to all odd integers are the unconditional result of deterministic time complexity $O(N^{1/3+\varepsilon})$ in [LN], and the heuristic result of deterministic time complexity $O(N^{1/4+\varepsilon})$ in [MK]. These results are derived from the continued fraction approximations of $q/p$ and $(q/p)^{1/2}$ respectively. In addition, the congruence congruence $2^{N+1} \equiv 2^{p+q} \bmod N$ is used in [P] and [MP] to obtain a $O(N^{1/3+\varepsilon})$ complexity for $N = pq$, $N^{1/3} < p < q$.

Another approach is explored here, it is simple, unconditional and rigorous. This approach utilizes the approximation

$$X_a = N^{1/2} + aN^{1/4} + \frac{N}{N^{1/2} + aN^{1/4}} \qquad (7)$$

of $p + q$. This approximation satisfies the inequality $| p + q - X_a | < 2N^{1/4}$ for any factor $N^{1/4} < p$, see (6).

**Theorem 4.** The *extended* Fermat method decomposes an odd integer $N$ in deterministic exponential time $O(N^{1/4+\varepsilon})$.

Proof: The aim is to compute a solution $x_0 = p + q$ of the equation $X^2 - Y^2 = 4N$. Without loss in generality let $N = pq$ where $p$ and $q$ are primes such that $N^{1/4} < p < q$, and write $p + q = X_a + t$, where $0 \le | a |, | t | < 2N^{1/4}$. Moreover, let $| T | < N^{1/4}$ be a primitive root modulo $N$. Now observe that for a product of two primes $N = pq$, the congruence $T^{N+1} \equiv T^{p+q}$ mod $N$ holds for any $T \in ( 1, N )$, $\gcd(T, N) = 1$. Applying the Shank's baby step giant step procedure to the congruence

$$T^{N+1} \equiv T^{X_a + t} \bmod N \qquad (8)$$

returns a solution $(a, t)$, see [MZ, p. 105] for implementation details. Clearly, this procedure has a deterministic the running time complexity of $O(N^{1/4+\varepsilon})$ arithmetic operations. ∎

A primitive root $| T | < O(N^{1/4+\varepsilon})$ modulo $N$ exists unconditionally, see Theorem 29, and $| T | < O((\log N)^6)$ conditioned on the Riemann hypothesis, see Theorem 31. In practice a significantly smaller one can be determined, see [BH], [PS], [B2]. In fact any element $T \in ( 1, N )$ of order $\mathrm{ord}_N(T) > 3N^{1/2}$ is sufficient. Also observe that there are other ways of using a solution $(a, t)$ to extract the factors of $N$. For example, by means of the equation $X^2 - x_0 X + N = 0$, or by means of $\gcd(T^{(N+1)/2} \pm T^{(p+q)/2}, N)$.

For a pair of factors $p$ and $q$ that share a common factor $d$, that is, $d \mid p - 1$ and $d \mid q - 1$ of size $d = O(N^\alpha)$, the complexity can be reduced to $O(N^{1/2-2\alpha})$ using standard technique of divisors in residue classes due to ? Landry 1869, Lehmer 1930 and other, see [W3].

A *sparse* integer $a = a_1 2^{v_1} + \cdots + a_k 2^{v_k}$ is a compact binary representation of some integers. These expansions are very effective in certain applications. The parameters $k$ and $v$ depend on the application.

**Theorem 5.** Let $N = pq$ and let $p = N^{1/2} + aN^{1/4} + b$. Then the following holds.
(i) If $a$ and $b$ are sparse integers, then the *extended* Fermat method has deterministic superpolynomial time complexity $O((\log N)^{O((\log \log N)^c)})$.
(ii) If $p \pm q = 2N^{1/2} + rN^{1/4} + s$ has sparse integer coefficients $r$ and $s$, then the *extended* Fermat method has deterministic superpolynomial time complexity.

Proof (i) Since $0 \leq |a|, |b| < N^{1/4}$, let $k = O(\log\log(N))^c$, and $0 \leq v_i \leq v = 1 + (\log N)/4$, $c > 0$ constant. Then this claim is verified by a direct calculation of an integral solution of

$$\frac{N}{N^{1/2} + aN^{1/4} + b} \tag{9}$$

as $a$ and $b$ vary over the set $A(k,v) = \{ \ a = a_1 2^{v_1} + \cdots + a_k 2^{v_k} \ \}$ of sparse integers parametized by $k$ and $v$. Similarly, the proof of (ii) uses a direct calculation of a solution $x_0 = p + q = 2N^{1/2} + rN^{1/4} + s$ or $y_0 = q - p$ of $X^2 - Y^2 = 4N$, such that $0 \leq |r|, |s| < N^{1/4}$. The complexity of the algorithm is dominated by the square of the cardinality $\#A(k,v) = (2v)^k \leq O((\log N)^{O((\log\log N)^c)})$. ∎

A lattice reduction theory version of this idea yields a result without any restriction on the coefficients $0 \leq |b|, |s| < N^{1/4}$, but at the expense of a more complicated algorithm.

The previous result is a very useful test to identify weak balanced integers $N = pq$, $p < q < 2p$. More precisely, a large class of balanced integers

$$\mathcal{E} = \{ \ N = pq : p = N^{1/2} + aN^{1/4} + b < q < 2p, \text{ or } p \pm q = 2N^{1/2} + rN^{1/4} + s \ \}, \tag{10}$$

where $a$, $b$, $r$, and $s$, are sparse integers decompose in deterministic superpolynomial time complexity. These integers are indistinguishable from random balanced integers, and can pass all the well known tests such as

(a) Have large sums and differences $p + q$ and $q - p \approx N^{1/2}$,
(b) The associated integers $p \pm 1$ and $q \pm 1$ are not smooth,
(c) The factors $p$ and $q$ have dense binary and nonadjacent form expansions,
(d) Are hard to factor using ordinary factoring algorithms, etc.

The example below demonstrates the concept of the Test/Factor, see example 8 too for a similar construction.

***Example* 6.** Test/factor the random balanced integer $N = pq = 448316072600119$, the 8 digits primes $p < q < 2p$ were selected at random.

Fix the parameters $k \leq \log\log(N)^2 = 5.60$, and $v_i \leq v = .25\log(N) = 12.16$. For each sparse integer $a = a_1 2^{v_1} + \cdots + a_k 2^{v_k}$ of weight $\text{wt}(a) = \#\{ \ a_i \neq 0 \ \} \leq k$, take the approximation

$$p_a = [N^{1/2}] + a[N^{1/4}] = 21173475 + 4601a$$

of $p \mid N$. A direct search for a solution of $X^2 - Y^2 = 4N$ is conducted using the approximation of $x_0 = p + q$ given by the integers

$$x = 21173475 + 4601a + N/(21173475 + 4601a) + t$$

as $t$ varies, $0 \leq |t| \leq \log(N)^2 = 1642$. A solution is detected by an integral value of $y = (x^2 - 4N)^{1/2}$. A few data points are shown in the table:

| $a_4 2^{v_4} + a_3 2^{v_3} + a_2 2^{v_2} + a_1 2^{v_1}$ | $= a$ | $t$ | $x = x_a + t$ | $Y$ |
|---|---|---|---|---|
| 0 | 0 | 1 | 42346952 | 7299.988219 |
| 1 | 1 | 3 | 42346955 | 17532.01497 |
| $2^2 + 1$ | 5 | 10 | 42346986 | 54156.1051 |
| … | … | … | | |
| $2^{11} - 2^8 - 2^6 - 2^2$ | 1724 | 339 | 44509024 | 13703610 |

The factors are $p = (44509024 - 13703610)/2 = 15402707$, and $q = (44509024 + 13703610)/2 = 29106317$. Accordingly, $N \in \mathcal{E}$ is in the class of balanced integers of deterministic superpolynomial time complexity with respect to $k$ and $v$.

## 3 Sparse Difference Methods

The *sparse difference method* decomposes integers whose prime differences (or sum) are sparse integers in superpolynomial time. This algorithm augments the Fermat factoring method. A sizable proportion of the integers factorable by these methods cannot be factored within reasonable time using the current best algorithms such as the NFS, etc., see the demonstration in example 8.

To initiate the development of this idea, consider the subset of integers with sparse differences:

$$\mathfrak{C} = \{ \, N = pq : q - p = a_1 2^{v_1} + \cdots + a_k 2^{v_k} \, , \, a_i \in \{-1, 0, 1\} \, \}, \tag{11}$$

where $p$, $q$ are primes and parameters $k$ and $v$ are fixed or a function of $N$. The associated counting function is

$$\mathfrak{C}(X) = \#\{ \, N = pq \leq X : q - p = a_1 2^{v_1} + \cdots + a_k 2^{v_k} \, , \, a_i \in \{-1, 0, 1\} \, \}. \tag{12}$$

The corresponding subset of sparse integers

$$A(X) = \{ \, a = a_1 2^{v_1} + \cdots + a_k 2^{v_k} : a_i = -1, 0, 1 \, \}, \tag{13}$$

is a complete set for the subset of balanced composites $\mathfrak{C}(X)$ in the sense that for every integer $N \in \mathfrak{C}(X)$, the prime difference $q - p \in A(X)$.

***Theorem 7.*** Let $c > 0$ be a constant, and let $X > 0$ be a large number. Then an integer $N \in \mathfrak{C}(X)$ can be factored in deterministic superpolynomial time $O((\log N)^{O((\log \log N)^c)})$.

Proof: Firstly observe that the ordered pair $(U_0 = q,\ V_0 = q - p)$ is a solution of the polynomial equation $U^2 - UV - N = 0$. Moreover, by hypothesis, $N \in \mathfrak{C}(X)$ implies that there exists a sparse integer $a \in A(X)$ such that the roots of the polynomial $U^2 - aU - N = 0$ are integral, in fact the factors of $N \leq X$.

The complexity of the algorithm is dominated by the calculations of the discriminant/roots of this equation, so the running time is dominated by $O(2 \log N)^k) \leq O((2 \log X)^k)$ arithmetic operations. $\blacksquare$

The asymptotic performance of the sparse difference method on the subset of integers $N = pq$ with $q - p = a_1 2^{v_1} + a_2 2^{v_2} + \cdots + a_k 2^{v_k}$ parametized by $k = O(\log\log(N)^c)$ and $v = O(\log N)$, $c > 0$ constant, surpasses the performance of any subexponential time factoring algorithms:

$$\left(2 \log N\right)^k \leq \left(\log N\right)^{c_1 (\log\log N)^c} = e^{c_1 (\log\log N)^{c+1}} < e^{c_0 (\log N)^\alpha (\log\log N)^{1-\alpha}}, \qquad (14)$$

where $0 < \alpha < 1$, $c$, $c_0$, $c_1 > 0$.

**Sparse Difference Algorithm**

Input: $N \in \mathfrak{C}(X)$.
Output: $p$, $q$.
1. Fix the parameters $k \leq c_0 (\log\log N)^c$, and $1 \leq v_i \leq c_1 \log(N)$, $c_0$, $c_1 > 0$ constants.
2. For $b = b_1 2^{u_1} + b_2 2^{u_2} + \cdots + b_k 2^{u_k}$, $b_i \in \{-1, 0, 1\}$ do
3. For $a = a_1 2^{v_1} + a_2 2^{v_2} + \cdots + a_k 2^{v_k} > 0$, $a_i \in \{-1, 0, 1\}$, compute $a^2 \pm 4bN$.
4. If $a^2 \pm 4bN$ is a square, then compute the root of $U^2 \pm aU \pm bN = 0$.
5. Return $-u_0 = p$, $u_1 = q$.

In step 3, the transform $p \pm q \rightarrow p \pm bq$ is included in the algorithm. This is motivated by simple cancellation property of binary addition/subtraction and it purpose is to explore a possible sparse integer solution $(R, S) = (a, b)$ of the equation $R = p \pm qS$. Quite often a nonsparse integer $p \pm q$ converts easily into a sparse integer by a simple shift $p \pm q \rightarrow p \pm 2^t q$, where $t \geq 1$.

In practice, the algorithm need to be optimized in order to remove the redundancy in the representations of the integers $a$. Ordering the digit patterns $a_k, \ldots, a_1 a_0$, by weight $\text{wt}(a) = \#\{ a_i \neq 0 \}$, $a_i \in \{ -1, 0, 1 \}$, $0 \leq i \leq k$, and leading digits $a_j = 1$ might accomplish this task.

***Example* 8.** The large balanced integer (251 digits) given below is easy to generate and passes all the well known tests: (i) have a large prime difference $q - p \approx N^{1/2}$, (ii) $p \pm 1$ and $q \pm 1$ are not smooth, (iii) it does not have any obvious patterns, (iv) it cannot be factored easily even with the best factoring algorithms. Nevertheless, *sparse difference* *algorithm* (Theorem 7) decomposes it in deterministic polynomial time, it requires less than $\log(N)^6$ arithmetic operations, try it.

$N = 48315390142927646144846003944586659635470343110930079906536801919808$
$\backslash 433501794456448193128039681227878838132098323453091986375810913119697 90$
$\backslash 5836154709673810134806193544846782452774974447854479640244139608220133 2$
$\backslash 17228623597586203471602828877526271389709.$

**Estimate of the Number of Certain Balanced Composites**
The effectiveness of the sparse difference method arises from the high density of primes pairs $p$ and $q$ with sparse differences or sums.

***Example* 9.** The balanced primes of $101 = p < q < 2p$ with sparse differences parametized by $k < \log\log 2p^2 = 3.84$, and $v_i < v = \log 2p^2 = 7.15$, are the following:

| | | |
|---|---|---|
| $101 + 2 = 103$ | $101 + 2^5 + 2 = 137$ | $101 + 2^6 + 2 = 167$ |
| $101 + 2^2 + 2 = 107$ | $101 + 2^5 + 2^2 + 2 = 139$ | $101 + 2^6 + 2^3 = 173$ |
| $101 + 2^3 = 109$ | $101 + 2^5 + 2^3 + 2 = 143$ | $101 + 2^6 + 2^4 - 2 = 179$ |
| $101 + 2^3 + 2 = 113$ | $101 + 2^5 + 2^4 = 149$ | $101 + 2^6 + 2^4 = 181$ |
| $101 + 2^4 + 2 = 119$ | $101 + 2^5 + 2^4 + 2 = 151$ | $101 + 2^7 - 2^5 + 1 = 193$ |
| $101 + 2^5 - 2^2 + 2 = 127$ | $101 + 2^5 + 2^2 + 2^3 = 157$ | $101 + 2^6 + 2^5 = 197$ |
| $101 + 2^5 - 2 = 131$ | $101 + 2^6 - 2 = 163$ | $101 + 2^6 + 2^5 + 2 = 199.$ |

The balanced primes of $103 = p < q < 2p$ includes distinct from those of 101, the balanced primes of $107 = p < q < 2p$ includes distinct from those of 101, and 103, etc.

The density of the subset of integers that have sparse prime differences should be asymptotic to something of the form

$$\mathfrak{C}(X) = \#\{\ N = pq \leq X : q - p = a_1 2^{v_1} + \cdots + a_k 2^{v_k}\ \} \sim \frac{X}{(\log X)^B}, \qquad (15)$$

for some constant $B \geq 2$, as $k \to \infty$, see Theorems 18, 19 and [GL] for supporting evidences. This is equivalent to proving that the subset of prime numbers in the subset of dePolignac numbers $q = p + a_1 2^{v_1} + \cdots + a_k 2^{v_k}$ has positive density.

Thus it is quite possible that the set $\mathfrak{C}(X)$ has nonnegligible density $C/(\log(X)^A, A > 0,$ relative to the set of balanced integers $\mathcal{B}(X) = \#\{\ N = pq \leq X : p < q < 2p\ \}$.

# 4 Sparse Exponent Method
The practicality of the proposed *sparse exponent method* is indirectly proportional to the largest prime factor of $p - 1$ or $q - 1$. This is indirectly quantified by the smoothness parameters $k$, $v$ and the random element $T$.

To establish the framework of the *sparse exponent method*, let $A(k,v) = \{ a = a_1 2^{v_1} + \cdots + a_k 2^{v_k} \}$ be the subset of sparse integers parametized by $k$ and $v$, and consider the system of inequations

$$(A_0p + B_0)(A_0p + B_1) \cdots (A_ip + B_{j-1})(A_ip + B_j) = c(q - 1), \tag{16}$$
$$(A_0p + B_0)(A_0p + B_1) \cdots (A_ip + B_{j-1})(A_ip + B_j) \neq d(p - 1),$$

where $A_i, B_j \in A(k,v)$ and $c, d \in \mathbb{Z}$. Observe that unlike the $p - 1$ method, the $p + 1$ method, and the elliptic curve method, the prime factors of the integers $A_ip + B_j$ are not restricted to a interval $[1, N^\beta]$, $0 < \beta < 1/4$, but can be significantly larger than $N^\beta$ depending on the parameters $k$ and $v \leq O(\log N)$. Thus there is a larger supply of primes $p_{i,j} \mid (A_ip + B_j)$ to construct a multiple of $p - 1$ or $p + 1$. The product of all the primes $p_{i,j}$ is accumulated in the effectively computable running exponent

$$E_{i,j} = (A_0N + B_0)(A_0N + B_1) \cdots (A_iN + B_{j-1})(A_iN + B_j), \tag{17}$$

where $0 \leq i, j \leq \#A(k,v)$. A more general *random exponent method* calls for a random running exponent

$$E_{i,j} = (A_0N + B_0)(A_0N + B_1) \cdots (A_iN + B_{j-1})(A_iN + B_j), \tag{18}$$

where $A_i, B_j \in \mathbb{Z}$ are random integers, $0 \leq i, j \leq M$, where $M$ is a fixed parameter.

***Theorem* 10.**   Suppose that the system (17) has a sparse integer solution, then $N = pq$ factors in probabilistic superpolynomial time.

Proof: Let $T$ be a random integer, $1 < \mid T \mid < N - 1$, and let $E_{i,j} = c(q - 1)$ and $E_{i,j} \neq d(p - 1)$ be satisfied by some $0 \leq i, j \leq \#A(k,v)$. Then it follows that

$$(A_0N + B_0)(A_0N + B_1) \cdots \quad (A_iN + B_{j-1})(A_iN + B_j) \equiv \tag{19}$$
$$(A_0p + B_0)(A_0p + B_1) \cdots \quad (A_ip + B_{j-1})(A_ip + B_j) \equiv 0 \bmod (q - 1)$$

and

$$(A_0N + B_0)(A_0N + B_1) \cdots \quad (A_iN + B_{j-1})(A_iN + B_j) \equiv$$
$$(A_0q + B_0)(A_0q + B_1) \cdots \quad (A_iq + B_{j-1})(A_iq + B_j) \not\equiv 0 \bmod (p - 1)$$

These imply that a nontrivial system of congruences

$$T^E \not\equiv 1 \bmod p, \text{ and } T^E \equiv 1 \bmod q,$$

where $E = E_{i,j}$ is the running exponent, can be generated in probabilistic superpolynomial time.  ∎

In almost all cases the function $f(E, N, T) = \gcd(T^E - 1, N)$ increases from 1 to $p$ or $q$ to $N$. It is rare to have $\gcd(T^E - 1, N) = N$. This later event occurs whenever the element $T \in \mathbb{Z}_N$ has (usually small) common order $r$ mod $p$ and mod $q$. In other word, the order $r = \mathrm{ord}_N(T)$ divides both integers $p - 1$ and $q - 1$.

### Sparse Exponent Algorithm

Input: $N = pq$.

Output: $p$, $q$.

1. Fix the parameters $k \leq c_0 \log\log(N)^c$, and let $0 \leq v_i \leq v = c_1 \log(N)$, where $c, c_0, c_1, > 0$ are constants.

2. For $0 \leq m \leq \log(N)$, select a random integer $T \in (1, N)$ such that $|T| > 1$.

3a. For $A_i, B_j \in A(k, v) = \{$ sparse integers parametized by $k, v \}$, $0 \leq i, j < \#A(k, v)$, do

3b. Compute $E_{i,j} = (A_i N + B_j) E_{i,j-1}$ the running exponent.

3c. Compute $T_{i,j} \equiv T_{i,j-1}^{E_{i,j}} \bmod N$, $T_{0,0} \equiv T^{E_{0,0}} \bmod N$.

3d. Compute $d_0 = \gcd(T_{i,j} - 1, N)$ and $d_1 = \gcd(T_{i,j} + 1, N)$.

4a. If $d_0$, or $d_1 = N$, then compute a square root $w$ of unity other than $\pm 1$.

4b. Compute $d_0 = \gcd(w - 1, N)$ and $d_1 = \gcd(w + 1, N)$.

5. If $1 < d_0 < N$ or $1 < d_1 < N$, then halt.

6. Return $d_i = p$, $q = N/p$.

**Remark:** In step 4, the condition $\gcd(a^E - 1, N) = N$ offers an opportunity to compute a square root of unity other than $\pm 1$ whenever $N$ is composite. The square root of unity modulo $N$ might be one of the integer in the series

$$w_0 \equiv a^s \bmod N, w_1 \equiv a^{2s} \bmod N, w_2 \equiv a^{2^2 s} \bmod N, ..., w_{r-1} \equiv a^{2^{r-1} s} \bmod N, \quad (20)$$

where $E = 2^r s$. There is also the possibility of computing an $n$th root of unity as an algebraic or primitive factor of $a^E - 1$ for some $n \mid E$.

***Example*** **11.** The examples below demonstrate a simple version of the sparse exponent method on sets of structured integers.

(i) A product $N = pq$ of Germain primes $p$, and $q = 2kp + 1$, with $k \geq 1$, and $2k(2k + 1) \not\equiv 0 \bmod (p - 1)$, decomposes in $k \leq O(\log(N)^c)$ steps. To verify this claim take $E = 2kN$ and $T = 2$. Then

$2kN \equiv 0 \bmod (q - 1)$ and $2kN \not\equiv 0 \bmod (p - 1)$.

Accordingly $2^{2kN} \equiv 1 \bmod q$ and $2^{2kN} \not\equiv 1 \bmod p \implies \gcd(2^{2kN} - 1, N) = 2kp + 1$.

(ii) A Mersenne number $N = 2^r - 1$ with $r$ prime. Put $N = pq = (2rk + 1)(2rm + 1)$, with $k < m$ and $km < 2^r/r^2$. Assume that the integers $k$ and $m$ satisfy the simple linear conditions

$Ak + B \equiv 0 \bmod m$ and $Am + B \not\equiv 0 \bmod k$. (21)

Now take $E = (N - 1)A + 2rB$ and $|T| > 2$. Then it follows that

$(N - 1)A + 2rB \equiv 2r(Ak + B) \equiv 0 \bmod (q - 1)$ and
$(N - 1)A + 2rB \equiv 2r(Am + B) \not\equiv 0 \bmod (p - 1)$.

Consequently, $T^{(N-1)A + 2rB} \equiv 1 \bmod q$, and $T^{(N-1)A + 2rB} \not\equiv 1 \bmod p$.

These data imply that $N > \gcd(T^{(N-1)A + 2rB} - 1, N) > 1$.

The number of steps needed to decompose $N$, depends on the number of sparse integers pairs $A$, $B$ required to realize the linear conditions (21).

(iii) A Fermat number $F_n = 2^{2^n} + 1$, $n \geq 0$. Put $N = pq = (2rk + 1)(2rm + 1)$, with $k < m$ and $km < 2^{2^n - 2n - 4}$. Assume that the integers $k$ and $m$ satisfy the simple linear conditions (21). Now take $E = (N - 1)A + 2^{n+2}B$ and $|T| > 2$. Then it follows that

$(N - 1)A + 2^{n+2}B \equiv 2^{n+2}(Ak + B) \equiv 0 \bmod (q - 1)$ and
$(N - 1)A + 2^{n+2}B \equiv 2^{n+2}(Am + B) \not\equiv 0 \bmod (p - 1)$.

Consequently, $T^{(N-1)A + 2^{n+2}B} \equiv 1 \bmod q$ and $T^{(N-1)A + 2^{n+2}B} \not\equiv 1 \bmod p$.

The number of steps needed to decompose $N$, depends on the number of sparse integers pairs $A$, $B$ required to realize the linear conditions (21).

For instance, $F_5 = 2^{32} + 1 = (2^7 \cdot 5 + 1)(2^7 \cdot 52347 + 1)$ decomposes in about three steps since $A = 1$, $B = -2$ is sufficient to satisfy (21). To factor it, choose a random $|T| > 2$, and use the expression $\gcd(T^{(N-1)A + 2^{n+2}B} \pm 1, N)$ with $A$, $B$ sparse variables.

## 5 Weak Balanced Integers and Hard to Factor Integers

A summary of the well known and (new) classes of weak integers are listed below. These integers are indistinguishable from random balanced integers, but have superpolynomial time complexities or other practical time complexities. These techniques are of interest in Cryptography and Standard X9.31. The construction of hard to factor balanced integers $N = pq$ requires a careful inspection of the following classes of weak integers.

Take $N = pq$, with $p$ and $q$ primes. The well known classes of weak integers are the followings.

(a) $p - 1, p + 1, q - 1, q + 1, p + 1 - a_p$, or $q + 1 - a_q$, where $| a_p | < 2p^{1/2}$, is $B$-smooth, where the parameter $B > 0$ is of practical size, for instance, $B \leq O(\log(N)^6)$ is practical today. The factorization of $N$ is handled by the $p \pm 1$ methods in $O(B \log(N)^c)$ arithmetic operations or the elliptic curve method in subexponential time.

(b) $N^{1/2} - N^{1/4} \leq p \leq N^{1/2}$ or $N^{1/2} \leq q. \leq N^{1/2} + N^{1/4}$.
The factorization of $N$ is handled by the standard Fermat method in deterministic polynomial time. Otherwise, $p < N^{1/2} - N^{1/4+\alpha}$, and $N^{1/2} + N^{1/4+\alpha} < q$ for some $\alpha > 0$, and the algorithm runs in $O(N^{2\alpha})$ arithmetic operations.

The (new) classes of weak balanced integers are the followings.

(c) $p = N^{1/2} + aN^{1/4} + b$,           where $| a | \leq O(N^\varepsilon)$ and $| b | \leq O(N^{1/4-\varepsilon})$, $\varepsilon > 0$.
The factorization of $N$ is handled by the extended Fermat method in $O(N^\varepsilon)$ arithmetic operations, see Theorem 3.

For instance, a 1024-bit integers $N = pq$, $p = N^{1/2} + aN^{1/4} + b < q < 2p$, for which $| a | \leq O(N^{1/20})$ and $| b | \leq O(N^{1/5})$ has a time complexity of $O(N^{1/20})$ arithmetic operations. This can be done effectively using current technology.

(d) $p = N^{1/2} + aN^{1/4} + b$,        where $a$ and $b$ are sparse integers, $0 \leq | a |, | b | < N^{1/4}$.
The factorization of $N$ is handled by the extended Fermat method in superpolynomial time, see Theorem 3 or 5.

(f) $p = 2N^{1/2} + rN^{1/4} + s$,        where $r$ and $s$ are sparse integers, $0 \leq r, s < N^{1/4}$.
The factorization of $N$ is handled by the extended Fermat method in superpolynomial time, see Theorem 3 or 5.

(g) $q = p + a_1 2^{v_1} + \cdots + a_k 2^{v_k}$,        where $a_i \in \{ -1, 0, 1 \}$, $0 \leq v_i < c \log(N)$.
The factorization of $N$ is handled by the *sparse difference method*, in superpolynomial time, see Theorem 7 and example 8.

To construct hard to factor integers and to prevent the pitfalls described before (weak integers $N = pq$), the expansions of the integers $a$, $b$, $r$, $s$, $p$, $q$, $p + q$, and $q - p$ must be dense. The $D$-expansions of these numbers must have large numbers of nonzero digits. Several $D$-expansions of these numbers using various digit sets $D = \{ 0, 1 \}$, $\{ -1, 0, 1 \}$, $\{ 0, 1, 3 \}$, etc., should be tested. However, it is likely that the two digits sets $\{ -1, 0, 1 \}$ and $\{ 0, 1, 3 \}$ are sufficient because these two digit sets minimize the weight of any expansion of an integer, see [SM]. Nevertheless, the selection of hard to factor integers could be more difficult than before.

## Sparse/Dense Expansions and Complexity

The sparse integer expansions are related to *non adjacent form* expansions of integers. The digits of a non adjacent form expansion $N = a_n 2^n + \cdots + a_1 2 + a_0$ satisfies $a_i a_{i+1} = 0$, and $a_i \in \{ -1, 0, 1 \}$, $0 \leq i \leq n$. The distribution of the number of nonzero digits $a_i = \pm 1$ is approximately Gaussian of mean $\mu = (1/3)\log_2(N)$, and standard deviation $\sigma = (22/108)\log_2(N)$, see [HP], and [SM] for the advanced theory.

The binary or nonadjacent form expansions of an integer $N$ has been described as *sparse* if the weight satisfies $\mathrm{wt}(N) = \#\{\ a_i \neq 0\ \} \leq O(\log\log(N)^c)$, $c > 0$. Otherwise, it has been called *dense* and $\mathrm{wt}(N) = \#\{\ a_i \neq 0\ \} \geq O(\log(N)^\alpha)$, $\alpha > 0$. This definition is inharmony with the standard measure of complexity $L(c, \alpha) = e^{c(\log N)^\alpha (\log\log N)^{1-\alpha}}$.

The application of the theory of integers expansions to integers factorization is an unexplored area of research. Properties of the expansions of integers (using single base or multibase etc.) could very well lead to an improvement on the time complexities of established results.

Roughly speaking, the nonadjacent form expansion of the integer $a = a_n 2^n + \cdots + a_1 2 + a_0 < N^{1/4}$ has an expected weight of $\mathrm{wt}(a) = \#\{\ a_i \neq 0\ \} \leq (\log(N))/12$. Thus the factorization of an arbitrary integer $N = pq$, $p = N^{1/2} + aN^{1/4} + b$, should have an expected time complexity of $O(N^{1/12})$ arithmetic operations. Furthermore, if the weight $\mathrm{wt}(a) = \#\{\ a_i \neq 0\ \} \leq (\log(N)^\alpha)$, $0 < \alpha < 1$, it has deterministic subexponential time complexity unconditionally. In contrast, the current best rigorous estimate has the time complexity $O(N^{1/5+\varepsilon})$ conditioned on the Riemann hypothesis, Shank 1974.

## 6 Notations And Background Materials

The bulk of integer factorization theory is focused on the class of difficult to factor set of balanced integers. The general purpose algorithms for factoring balanced integers are complex and run in subexponential times to exponential times. On the other hand, the algorithms for factoring nonbalanced large integers are simpler, more effective and have lower complexity.

**Balanced Integers**

In the investigation of the complexity of integer factoring it is sufficient to consider the running time complexity $\theta(N)$ of decomposing an arbitrary integer $N$ into two large factors. Moreover, since an arbitrary integer $N$ has fewer than $2\log(N)$ divisors the overall running time complexity of a complete factorization is $O(\theta(N)\log(N)^A)$, $A > 0$.

If the smallest prime factor $p$ of $N$ satisfies $p > N^{1/4}$, then it immediately follows that $N = p$ or $N = pq$ or $N = pqr$, so it has at most three prime factors. Similarly, if the smallest prime factor $p$ of $N$ satisfies $p > N^{1/3}$, then it immediately follows that $N = p$ or $N = pq$, so it has at most two prime factors. The most important and difficult case in integer factorizations is the class of 2-balanced integers $N = pq$ whose prime factors satisfy

$$p < N^{1/2} - N^{1/4+\alpha} \ \text{ and } \ N^{1/2} + N^{1/4+\alpha} < q < 2p, \text{ where } \alpha > 0.$$

***Definition* 12.** Let $a > 1$ be a small number, and let $p$ and $q$ be primes. The composite integer $N = pq$ is called *a*-balanced if the primes factors satisfy $p < q < ap$.

***Proposition* 13.** Let $N$ be an $a$-balanced integer, then $\sqrt{N} + N^{\delta} < q < \sqrt{aN}$, some $\delta \geq 0$.

Proof: Let $q$ be the larger factor of $N$. Then $N/q < \sqrt{N} < q < aN/q$. Multiplying across the inequalities by $q$ confirms the claim. ∎

***Proposition* 14.** If the integer $N = pq$ is $a$-balanced then $2\sqrt{N} < p + q < (1 + 1/a)\sqrt{aN}$.

Proof: Multiplying the inequalities $p < q < ap$ by $p$ yields $\sqrt{N/a} < p < \sqrt{N}$. Hence the image of the function $f(x) = x + N/x$ over the interval $[\sqrt{N/a}, \sqrt{N}]$ is precisely $[2\sqrt{N}, (1 + 1/a)\sqrt{aN}]$. ∎

***Proposition* 15.** Let $a, b > 1$ be numbers. If the primes factors of the integer $N = pq$ satisfy $q/a < p < q < bp$, then $(\sqrt{b} - 1)\sqrt{N} < q - p < (\sqrt{b} - \sqrt{1/a})\sqrt{N}$.

Proof: Similar algebraic manipulations as above. ∎

These inequalities provide a very sharp lower and upper estimates of $p$ and $q$. Namely, $\sqrt{N/a} < p < \sqrt{N}$ in Proposition 14, and combining Propositions 14 and 15 yield $.5(1 + \sqrt{2})\sqrt{N} < q < \sqrt{2N}$.

**Smooth Integers**

The set of $z$-smooth numbers is defined by $\Psi(z) = \{\, n \in \mathbb{N} : n \equiv 0 \bmod p \implies p < z \,\}$, the corresponding counting function is defined by

$\Psi(x, z) = \#\{\, n \leq x : n \equiv 0 \bmod p \implies p < z \,\}$.

The complement set of $z$-nonsmooth numbers is $\Phi(z) = \{\, n \leq x : n \not\equiv 0 \bmod p \impliedby p < z \,\}$ and its corresponding counting function is defined by

$\Phi(x, z) = \#\{\, n \leq x : n \not\equiv 0 \bmod p \impliedby p < z \,\}$.

***Theorem* 16.** Let $0 < z < x$ be real numbers. Then $\Psi(x, z) = cx^{\delta}(\log z)e^{-\log x/\log z}$ as $x \to \infty$ and $\delta > 1/2$.

Proof: As given in [CU, p. 68]. Start from definition and expand it:

$$\Psi(x, z) = \sum_{\substack{n \leq x \\ p \mid n \Rightarrow p < z}} 1 \ \leq \sum_{\substack{n \leq x \\ p \mid n \Rightarrow p < z}} \left(\frac{x}{n}\right)^{\delta} \leq x^{\delta} \prod_{p < z} (1 - 1/p^{\delta})^{-1}, \tag{22}$$

where $\delta$ is a small number. Rewrite the product as

$$\prod_{p<z}(1-1/p^\delta)^{-1} = c_0 \prod_{p<z}(1+1/p^\delta)\prod_{p<z}(1-1/p^{2\delta})^{-1} = c_1 \prod_{p<z}(1+1/p^\delta)$$

which is convergent for all $\delta > 1/2$. Thus $\Psi(x,z) = c_1 x^\delta \prod_{p<z}(1+1/p^\delta)$. A second elementary step converts it to the stated form. ∎

Often the lower estimate $\Psi(x,z) > x^{1-\log\log x/\log z}$ is adequate.

**Theorem 17.** Let $0 < z < x$ be real numbers, then
$$\Phi(x,z) = x^\delta \prod_{p<z}(1-1/p) + O(x(\log z)^2 e^{-\log x/\log z}) \text{ as } x \to \infty.$$

**Time Complexities of Algorithms**

An algorithm runs in *polynomial time* if there is a constant $c > 0$ such that $(\log N)^c$ is an upper estimate on the number of arithmetic operations required by the algorithm. This estimate is based on the worst case scenario.

The function $L(c,\alpha) = e^{c(\log N)^\alpha(\log\log N)^{1-\alpha}}$ interpolate between *polynomial* time at $\alpha = 0$, and *subexponential* time at $0 < \alpha < 1$, to *exponential* time at $\alpha = 1$.

A subexponential integer factoring algorithm has a run time complexity of the form

$$L(c_0,\alpha) = e^{c_0(\log N)^\alpha(\log\log N)^{1-\alpha}}, \tag{23}$$

where $0 < \alpha < 1$, and $c_0 > 0$ are constant. In comparison, a *superpolynomial* time integer factoring algorithm has a run time complexity of the form

$$M(b,c_1) \leq (\log N)^{c_1(\log\log N)^c} = e^{c_1(\log\log N)^{c+1}}, \ b, \ c_1, \ c > 0 \text{ constants.} \tag{24}$$

The current integer factoring algorithms are rated as follows.

(1) Trial divisions, the oldest integer factoring method, this is deterministic and rigorously analyzed. There are various versions, but the best time complexity achievable is about $O(N^{1/2}/\log(N)^A)$ arithmetic operations, some constant $A > 0$.

(2) The Continued Fraction Factoring Algorithm has a heuristic random time complexity of $O(e^{\sqrt{2(\log N)(\log\log N)}})$ arithmetic operations on an arbitrary odd integer $N$.

(3) The Quadratic Sieve Factoring Algorithm has a heuristic random time complexity of $O(e^{\sqrt{(\log N)(\log\log N)}})$ arithmetic operations on an arbitrary odd integer $N$.

(4)    The Number Fields Sieve Factoring Algorithm has a heuristic random time complexity of $O(e^{c(\log N)^{1/3}(\log\log N)^{2/3}})$ operations on an arbitrary odd integer $N$, where $0 < c \leq 1.923$.

(5)    The Elliptic curve method has a heuristic random time complexity of $O(e^{c(\log N)^{1/2}(\log\log N)^{1/2}})$ arithmetic operations on an arbitrary odd integer $N$, where $0 < c \leq 1.923$.

(6)    The $p - 1$ method, $p + 1$ method have probabilistic time complexity of about $O(N^{1/2}/\log(N)^4)$ arithmetic operations on an arbitrary odd integer $N$.

**Universal Exponent Method**

The universal exponent method is a general factoring technique in finite rings that support some form of unique factorizations and a gcd algorithm (it probably works in other cases too). This technique is derived from the orders of the elements in the multiplicative group of the finite rings. Legendre's Theorem states that the order of an element in a finite group is a divisor of the order of the group. The basic idea of the universal exponent method seems to be of unknown origin, see [S], but it can be traced back to at least the early 1900, see [W1]. The best known instances of this technique are the $p - 1$ method, the $p + 1$ method, (more generally the $\Phi_k(p)$ method) and the elliptic curve method.

The running time complexities of these algorithms are tied up to the sizes of the prime factors $p$ of the integer $N$ under test and the structures of the corresponding finite groups. Specifically:

(i) The effectiveness of the the $p - 1$ method is directly proportional to the largest prime factor of $p - 1$. This is directly quantifies by the smoothness parameter $B$ and the random element $a$.

(ii) The effectiveness of the $p + 1$ method is directly proportional to the largest prime factor of $p + 1$. This is directly quantified by the smoothness parameter $B$ and the random element $a$.

(iii) The effectiveness of the elliptic curve method is indirectly proportional to the largest prime factor of $p + 1 - a_p$ or $q + 1 - a_q$, where $|a_p| < 2p^{1/2}$. This is indirectly quantified by the smoothness parameter $B$, and the random elliptic curve $E$.

(iv) The effectiveness of the *sparse exponent method* is indirectly proportional to the largest prime factor of $p - 1$ or $q - 1$. This is directly quantified by the smoothness parameter $k$ and $v$ and the random element $a$.

**Classification of Algorithms**

An integer factorization algorithm will be classified as either multiplicative or additive. A *multiplicative* algorithm produces the factors of on integer by a process of multiplications, and an *additive* algorithm produces the factors of on integer by a process of additions.

(i) The Fermat factoring method, circa 1665, is an additive algorithm. It constructs a solution of the equation $X^2 - Y^2 = N$ by completing the squares through addition or equivalently by sequentially incrementing the initial approximation $X = \sqrt{N}$ of $(p + q)/2$.

(ii) The Seelhoff factoring method, see [W2, p. 128], and its modern versions such as the continued fraction method, random squares method, and number field sieve, et cetera, are multiplicative algorithms. These algorithms generate the solutions of the congruence $X^2 - Y^2 \equiv 0 \bmod N$ by completing the squares through multiplications.

(iii) The *sparse exponent method* is a combination additive/multiplicative algorithm. It determines a nontrivial square root of unity by a process of addition/multiplication. This algorithm augments the $p - 1$ method, the $p + 1$ method, and the elliptic curve method.

## 7 Theoretical Foundations

This section is a pointer to mathematical foundations of the new algorithm. The basic underpinning of the sparse difference method and the sparse exponent method rest on the theory of integer representations such as sums primes, powers, and composites. Much of the complicated theory of integer representations of has been known for many years, and emerged in response to challenges of the twin primes problem, dePolignac problem (1849), and the likes.

**Sum of Primes**

The dePolignac problem states that every odd integer is of the form $n = p + 2^v$. This is one of the simplest sum of primes problem. The density of these numbers has been settled by a sieve argument and other means.

Let $x > 0$, and $r(n)$ be the number of representation of an integer $n \leq x$ as $n = p + a^v$, where $|a| > 1$. The following sums

$$c_0 x \leq \sum_{n \leq x} r(n) \leq c_1 x \quad \text{and} \quad \sum_{n \leq x} r(n)^2 = c_2 x \tag{25}$$

are used in a density proof of dePolignac problem. The proof of the first sum in (24) is quite simple, a few lines, but the proof of the second sum is considerably longer, about two pages, see [NT, p. 199-204].

***Theorem 18.*** ([RV]) The set of integers $R(x) = \#\{ p + 2^v \leq x : p \text{ prime and } v \geq 0 \}$ has a cardinality of $R(x) = cx$, for some constant $c > 0$.

Proof: As given in [NT, p. 204]. Use the Cauchy-Schwarz inequality to get

$$(c_0 x)^2 \leq \left( \sum_{n \leq x} r(n) \right)^2 \leq R(x) \sum_{n \leq x} r(n)^2 \leq c_2 x R(x). \tag{26}$$

Thus for sufficiently large $x$, $R(x) \geq cx$, where $c > 0$ is a constant. ∎

Recently a few authors have given an explicit range, namely, $.1866x \leq R(x) \leq .9819x$, see [HR].

The Romanoff set $R(x)$ contains a surprisingly large number of primes, even for small $x > 0$. The subset $P(x) \subset R(x)$ of primes appear to have a cardinality of the form

$$P(x) = \#\{\, p + 2^v \leq x : p \text{ prime and } p + 2^v \text{ are primes}, v \geq 0 \,\} = rx/(\log x)^2, \tag{27}$$

with $r > 0$ constant. This observation seems to follows from the next result.

***Theorem*** 19.   Let $d \geq 2$ be a fixed even integer. The number of primes in the subset of primes $\{\, p \leq x : p \text{ and } | \, p + d \, | \text{ are primes} \,\}$ is given by

$$B(x) = \frac{cx}{(\log x)^2} \prod_{p \mid d} (1 + 1/p), \, c > 0 \text{ constant.} \tag{28}$$

The proof is derived from the Brun's sieve or other sieve, see [CU, p. 102] and the literature.

The Fermat factoring method is associated with the sequence of prime pairs $p$ and $p + 2n$ with $2n \leq X^{1/4}$. To determine the size $\mathcal{F}(X)$ of the domain of this method, it is essential to measure the magnitude of the set $B(X^{1/2}) = \#\{\, p \leq X^{1/2} : p \text{ and } p + 2n \text{ are primes, and } 2n \leq X^{1/4} \,\}$

Since the asymptotic expression

$$\sum_{n=1}^{X^{1/4}/2} \prod_{p \mid 2n} (1 + 1/p) \sim \sum_{n=1}^{X^{1/4}/2} \log(2n) \sim X^{1/4}, \tag{29}$$

holds, it follows that

$$\mathcal{F}(X) \leq \sum_{n=1}^{X^{1/4}/2} \frac{cX^{1/2}}{(\log X^{1/2})^2} \prod_{p \mid 2n} (1 + 1/p) \sim \frac{cX^{3/4}}{(\log X^{1/2})^2}. \tag{30}$$

**Quasiharmonic Sums**

***Theorem* 20.**   For $x \geq 2$, the following hold.

(i)   $\displaystyle \sum_{k \leq x,\, \gcd(k,n)=1} \frac{1}{k} = \frac{\varphi(n)}{n} \log x + C$ , where $C$ is a constant.

(ii)   $\displaystyle \sum_{k \leq x,\, k \equiv a \bmod n} \frac{1}{k} = \frac{\log x}{\varphi(n)} + C + O(1/\log x)$ , where $C$ is a constant.

An elementary proof of the first part appears in [RM, p.193].

***Theorem* 21.**   (Mertens Formulae 1874)   For $x \geq 2$, the following hold.

(i)   $\displaystyle \sum_{p \leq x,\, p \equiv a \bmod n} \frac{1}{p} = \frac{\log \log x}{\varphi(n)} + c(a,n) + O(1/\log x)$ , where $c(a, n)$ is a constant.

(ii)   $\displaystyle \prod_{p \leq x,\, p \equiv a \bmod n} (1 - 1/p) = \frac{e^{-\gamma}}{\varphi(n) \log x} + O(1/(\log x)^2)$ , where $\gamma$ is Euler's constant.

See [NZ, p. 128].

**Distribution of Primes Numbers**

A prime number has two divisors 1 and itself.

***Theorem* 22.**   (Fundamental Theorem of Arithmetic)   Every integer has a unique factorization of the form $N = p_1^{e_1} p_2^{e_2} \cdots p_\omega^{e_\omega}$ , where $p_i$ is prime, $e_i \geq 1$, and $\omega = \omega(N) \geq 1$.

***Theorem* 23.**   (Euclid, 300BC)   There are infinitely many primes.

***Theorem* 24.**   (Ramanujan 1936)   Let $\omega(N) = \#\{$ prime $p : p \mid N \}$, and let $\varepsilon > 0$. Then $\left| \omega(N) - \log \log N \right| \leq \left( \log \log x \right)^{1/2 + \varepsilon}$ for almost all integers $N \leq x$.

***Theorem* 25.**   (Prime Number Theorem)   Let $x \geq x_0 > 0$. Then

(i)   $\dfrac{x}{2 \log x} < \pi(x) < \dfrac{2x}{\log x}$ ,

(ii)   $\pi(x) = li(x) + O(x e^{-c\sqrt{\log x}})$ , $c > 0$ constant,

(iii)   The probability of a random integer $n \leq x$ of being prime is approximately

$\dfrac{1}{2 \log x} < \dfrac{\pi(x)}{x} < \dfrac{2}{\log x}$ .

Statement (i) is due to Chebyshev 1850. The prime number theorem, statement (ii) was independently proved by delaVallee Poussin and Hadamard in 1886. The logarithmic integral approximation

$$\pi(x) \approx li(x) = \int_2^x \frac{dt}{\ln t} = x\left(\frac{1}{\ln x} + \frac{1}{(\ln x)^2} + \frac{1}{(\ln x)^3} + \cdots\right), \quad (31)$$

where ln is the logarithmus naturalis, was introduced by Gauss in 1792. It gives a probabilistic estimate of the prime counting function $\pi(x)$ with probability density function $1/\ln(x)$.

Conditioned on the Riemann hypothesis, the error terms of various results drop to square root order of magnitude. For example, the Prime Number Theorem is rephrased as

$$\pi(x) = li(x) + O(x^{1/2+\varepsilon}). \quad (32)$$

The Dirichlet result for the number of primes in arithmetic progression was subsequently refined.

***Theorem 26.*** (delaVallee Poussin 1896)   The number of primes $p \leq x$ in the arithmetic progression $nz + a$ is given by

$$\pi(x,n,a) = \frac{1}{\varphi(n)} li(x) + O(xe^{-c(\log x)^{3/6}(\log\log x)^{-1/5}}), \quad (33)$$

where $c > 0$ is a constant.

The primes are uniformly distributed in each equivalent class mod $n$ for $n \leq O(\log x)^A$, $A > 0$, consult Walfiz 1936.

***Theorem 27.*** (Linnik 1944)   The least prime $p(n)$ in the arithmetic progression $nz + a$ satisfies $p(n) \leq O(n^L)$, with $L > 0$ constant.

Currently, it is known that the value $L \leq 11/2$, consult Heath-Brown 1992.

***Conjecture 28.*** (Generalized Prime Number Theorem)   Let $r_1, r_2, ..., r_k$ be an admissible set of integers. Then
(i) There are infinitely many $k$-tuples of primes $n + r_1, n + r_2, ..., n + r_k$ as $n$ ranges over the integers $\geq 0$.
(ii) The number of primes is given by $\pi(x, r_1, ..., r_k) \sim c(R)x/(\log x)^k$, where $c(R) = \prod_{p \geq 2}(1 - v(R)/p)(1/1/p)^k$.

This is better known as Hardy-Littlewood conjecture. Here the constant $c(R) > 0$ whenever the set $R = \{r_1, r_2, ..., r_k\} \not\equiv \mathbb{Z}_p$ for all prime $p$. Otherwise $c(R) = 0$.

**Primitive Roots**

***Theorem 29.*** ([BR])  For each $\varepsilon > 0$ and for every prime $p$, every interval of length $H > p^{1/4+\varepsilon}$ contains $\dfrac{\varphi(p-1)}{p-1} H (1 + O(p^{-\delta}))$ primitive roots mod $p$, where $\delta > 0$ depends only on $\varepsilon$. In particular, it follows that the least primitive root $g$ mod $p$ satisfies the estimate $g = O(p^{1/4+\varepsilon})$.

The average order of magnitude of the least primitive root $g(p)$ modulo $p$ is as follows.

***Theorem 30.*** ([BE])  For large $X$, $\pi(x)^{-1} \sum_{p \leq x} g(p) \prec\prec (\log x)^2 (\log \log x)^4$ the summation being extended over prime numbers $p$.

***Theorem 31.*** ([SP])  Assuming the Riemann Hypothesis, the least primitive root mod $p$ is $g(p) = O(r^4 (1 + \log r)^4 \log(p)^2)$, where $r = \omega(p-1)$.